\documentclass[preprint,12pt]{elsarticle}

\usepackage{amssymb}
\usepackage{amsthm}
\usepackage{graphics}
\usepackage{epsfig}
\usepackage{graphicx}
\usepackage{subfigure}
\usepackage{tikz}
\usetikzlibrary{calc}
\newcommand{\factorial}[1]{%
\tikz[baseline]
{\node[anchor=base,inner sep=0.3ex](mynode)
{\ensuremath{#1}};\draw(mynode.north west)--(mynode.south west)--
(mynode.south east);\path[use as bounding box]($(mynode.south west)+(-0.3ex,-0.3ex)$)
rectangle($(mynode.north east)+(0.3ex,0.3ex)$);}
}
\usepackage{a4wide}
\usepackage{amsmath}
\usepackage{amsfonts}
\usepackage{enumerate}
\newcommand{\pf}{\noindent {\bf Proof: }}

\newtheorem{lemma}{\bf Lemma}[section]
 \newtheorem{example}{\bf Example}[section]
\newtheorem{theorem}{\bf Theorem}[section]

\newtheorem{corollary}[lemma]{\bf Corollary}

\newtheorem{remark}{\bf Remark}[section]

\journal{Communications in Algebra}

\begin{document}

\begin{frontmatter}



\title{Non-Zero Component Graph of a Finite Dimensional Vector Space\footnote{Dedicated to Professor Mridul Kanti Sen}}



\author{Angsuman Das\corref{cor1}}
\ead{angsumandas@sxccal.edu}

\address{Department of Mathematics,\\ St.Xavier's College, Kolkata, India.\\angsumandas@sxccal.edu}
\cortext[cor1]{Corresponding author}


\begin{abstract}
In this paper, we introduce a graph structure, called non-zero component graph on finite dimensional vector spaces. We show that the graph is connected and find its domination number and independence number. We also study the inter-relationship between vector space isomorphisms and graph isomorphisms and it is shown that two graphs are isomorphic if and only if the corresponding vector spaces are so. Finally, we determine the degree of each vertex in case the base field is finite.
\end{abstract}


\end{frontmatter}



\section{Introduction}
The study of graph theory, apart from its combinatorial implications, also lends to characterization of various algebraic structures. The benefit of studying these graphs is that one may find some results about the algebraic structures and vice versa. There are three major problems in this area: (1) characterization of the resulting graphs, (2) characterization of the algebraic structures with isomorphic graphs, and (3) realization of the connections between the structures and the corresponding graphs. 

The first instance of such work is due to Beck \cite{beck} who introduced the idea of zero divisor graph of a commutative ring with unity. Though his key goal was to address the issue of colouring, this initiated the formal study of exposing the relationship between algebra and graph theory and at advancing applications of one to the other. Till then, a lot of research, e.g., \cite{survey2,zero-divisor-survey,anderson-livingston,graph-ideal,power1,power2,mks-ideal,badawi} has been done in connecting graph structures to various algebraic objects. Recently, intersection graphs associated with subspaces of vector spaces were studied in \cite{int-vecsp-2,int-vecsp-1}. However, as those works were a follow up of intersection graphs, the main linear algebraic flavour of characterizing the graph was missing. 

Throughout this paper, vector spaces are finite dimensional over a field  $\mathbb{F}$ and $n=dim_{\mathbb{F}}(\mathbb{V})$. In this paper, we define a graph structure on a finite dimensional vector space $\mathbb{V}$ over a field $\mathbb{F}$, called Non-Zero Component Graph of $\mathbb{V}$ with respect to a basis $\{\alpha_1,\alpha_2,\ldots,\alpha_n\}$ of $\mathbb{V}$, and study the algebraic characterization of isomorphic graphs and other related concepts.

\section{Definitions and Preliminaries}
In this section, for convenience of the reader and also for later use, we recall some definitions, notations and results concerning elementary graph theory. For undefined terms and concepts the reader is referred to \cite{west-graph-book}.

By a graph $G=(V,E)$, we mean a non-empty set $V$ and a symmetric binary relation (possibly empty) $E$ on $V$. The set $V$ is called the set of vertices and $E$ is called the set of edges of $G$. Two element $u$ and $v$ in $V$ are said to be adjacent if $(u,v) \in E$. $H=(W,F)$ is called a {\it subgraph} of $G$ if $H$ itself is a graph and $\phi \neq W \subseteq V$ and $F \subseteq E$. If $V$ is finite, the graph $G$ is said to be finite, otherwise it is infinite. The open neighbourhood of a vertex $v$, denoted by $N(v)$, is the set of all vertices adjacent to $v$. A subset $I$ of $V$ is said to be {\it independent} if any two vertices in that subset are pairwise non-adjacent. The {\it independence number} of a graph is the maximum size of an independent set of vertices in $G$. A subset $D$ of $V$ is said to be {\it dominating set} if any vector in $V \setminus D$ is adjacent to at least one vertex in $D$. The {\it dominating number} of $G$, denoted by $\gamma(G)$ is the minimum size of a dominating set in $G$. A subset $D$ of $V$ is said to be a {\it minimal dominating set} if $D$ is a dominating set and no proper subset of $D$ is a dominating set. Two graphs $G=(V,E)$ and $G'=(V',E')$ are said to be {\it isomorphic} if $\exists$ a bijection $\phi: V \rightarrow V'$ such that $(u,v) \in E \mbox{ iff } (\phi(u),\phi(v)) \in E'$. A {\it path} of length $k$ in a graph is an alternating sequence of vertices and edges, $v_0,e_0,v_1,e_1,v_2,\ldots, v_{k-1},e_{k-1},v_k$, where $v_i$'s are distinct (except possibly the first and last vertices) and $e_i$ is the edge joining $v_i$ and $v_{i+1}$. We call this a path joining $v_0$ and $v_{k}$. A graph is {\it connected} if for any pair of vertices $u,v \in V,$ there exists a path joining $u$ and $v$.  The {\it distance} between two vertices $u,v \in V,~ d(u,v)$ is defined as the length of the shortest path joining $u$ and $v$, if it exists. Otherwise, $d(u,v)$ is defined as $\infty$. The {\it diameter} of a graph is defined as $diam(G)=\max_{u,v \in V}~ d(u,v)$, the largest distance between pairs of vertices of the graph, if it exists. Otherwise, $diam(G)$ is defined as $\infty$.

\section{Non-Zero Component Graph of a Vector Space}
Let $\mathbb{V}$ be a vector space over a field $\mathbb{F}$ with $\{\alpha_1,\alpha_2,\ldots,\alpha_n\}$ as a basis and $\theta$ as the null vector. Then any vector $\mathbf{a} \in \mathbb{V}$ can be expressed uniquely as a linear combination of the form $\mathbf{a}=a_1\alpha_1+a_2\alpha_2+\cdots+a_n\alpha_n$. We denote this representation of $\mathbf{a}$ as its basic representation w.r.t. $\{\alpha_1,\alpha_2,\ldots,\alpha_n\}$. We define a graph $\Gamma(\mathbb{V}_\alpha)=(V,E)$ (or simply $\Gamma(\mathbb{V})$) with respect to $\{\alpha_1,\alpha_2,\ldots,\alpha_n\}$ as follows: $V=\mathbb{V}\setminus \{\theta\}$ and for $\mathbf{a},\mathbf{b} \in V$, $\mathbf{a} \sim \mathbf{b}$ or $(\mathbf{a},\mathbf{b}) \in E$ if $\mathbf{a}$ and $\mathbf{b}$ shares at least one $\alpha_i$ with non-zero coefficient in their basic representation, i.e., there exists at least one $\alpha_i$ along which both $\mathbf{a}$ and $\mathbf{b}$ have non-zero components. Unless otherwise mentioned, we take the basis on which the graph is constructed as $\{\alpha_1,\alpha_2,\ldots,\alpha_n\}$.

{\theorem \label{basis-independent-theorem} Let $\mathbb{V}$ be a vector space over a field $\mathbb{F}$. Let $\Gamma(\mathbb{V}_\alpha)$ and $\Gamma(\mathbb{V}_\beta)$ be the graphs associated with $\mathbb{V}$ w.r.t two bases $\{\alpha_1,\alpha_2,\ldots,\alpha_n\}$ and $\{\beta_1,\beta_2,\ldots,\beta_n\}$ of $\mathbb{V}$. Then $\Gamma(\mathbb{V}_\alpha)$ and $\Gamma(\mathbb{V}_\beta)$ are graph isomorphic.}\\
\pf Since, $\{\alpha_1,\alpha_2,\ldots,\alpha_n\}$ and $\{\beta_1,\beta_2,\ldots,\beta_n\}$ are two bases of $\mathbb{V}$, $\exists$ a vector space isomorphism $T: \mathbb{V} \rightarrow \mathbb{V}$ such that $T(\alpha_i)=\beta_i, \forall i \in \{1,2,\ldots,n\}$.

We show that the restriction of $T$ on non-null vectors of $\mathbb{V}$, $\mathbf{T}:\Gamma(\mathbb{V}_\alpha) \rightarrow \Gamma(\mathbb{V}_\beta)$ is a graph isomorphism. Clearly, $\mathbf{T}$ is a bijection. Now, let $\mathbf{a}=a_1\alpha_1+a_2\alpha_2+\cdots+a_n\alpha_n; \mathbf{b}=b_1\alpha_1+b_2\alpha_2+\cdots+b_n\alpha_n$ with $\mathbf{a}\sim \mathbf{b}$ in $\Gamma(\mathbb{V}_\alpha)$. Then, $\exists~ i \in \{1,2,\ldots,n\}$ such that $a_i\neq 0,b_i\neq 0$. Also, $\mathbf{T}(\mathbf{a})=a_1\beta_1+a_2\beta_2+\cdots+a_n\beta_n$ and $\mathbf{T}(\mathbf{b})=b_1\beta_1+b_2\beta_2+\cdots+b_n\beta_n$. Since, $a_i\neq 0,b_i\neq 0$, therefore 
$\mathbf{T}(\mathbf{a}) \sim \mathbf{T}(\mathbf{b})$ in $\Gamma(\mathbb{V}_\beta)$.

Similarly, it can be shown that if $\mathbf{a}$ and $\mathbf{b}$ are not adjacent in $\Gamma(\mathbb{V}_\alpha)$, then $\mathbf{T}(\mathbf{a})$ and $\mathbf{T}(\mathbf{b})$ are not adjacent in $\Gamma(\mathbb{V}_\beta)$. \qed

{\remark The above theorem shows that the graph properties associated of $\Gamma$ does not depend on the choice of the basis $\{\alpha_1,\alpha_2,\ldots,\alpha_n\}$. However, two vertices may be adjacent with respect to one basis but non-adjacent to some other basis as shown in the following example: Let $\mathbb{V}=\mathbb{R}^2,\mathbb{F}=\mathbb{R}$ with two bases $\{\alpha_1=(1,0),\alpha_2=(0,1)\}$ and $\{\beta_1=(1,1),\beta_2=(-1,1)\}$.  Consider $\mathbf{a}=(1,1)$ and $\mathbf{b}=(-1,1)$. Clearly $\mathbf{a}\sim \mathbf{b}$ in $\Gamma(\mathbb{V}_\alpha)$, but $\mathbf{a}\not\sim \mathbf{b}$ in $\Gamma(\mathbb{V}_\beta)$.}

\section{Basic Properties of $\Gamma(\mathbb{V})$}
In this section, we investigate some of the basic properties like connectedness, completeness, independence number, domination number of $\Gamma(\mathbb{V})$.
{\theorem $\Gamma(\mathbb{V}_\alpha)$ is connected and $diam(\Gamma)=2$.}\\
\pf Let $\mathbf{a},\mathbf{b} \in V$. If $\mathbf{a}$ and $\mathbf{b}$ are adjacent, then $d(\mathbf{a},\mathbf{b})=1$. If $\mathbf{a}$ and $\mathbf{b}$ are not adjacent, since $\mathbf{a},\mathbf{b} \neq \theta$, $\exists \alpha_i, \alpha_j$ which have non-zero coefficient in the basic representation of $\mathbf{a}$ and $\mathbf{b}$ respectively. Moreover, as $\mathbf{a}$ and $\mathbf{b}$ are not adjacent, $\alpha_i \neq \alpha_j$. Consider $\mathbf{c}=\alpha_i + \alpha_j$. Then, $\mathbf{a}\sim \mathbf{c}$ and $\mathbf{b} \sim \mathbf{c}$ and hence $d(\mathbf{a},\mathbf{b})=2$. Thus, $\Gamma$ is connected and $diam(\Gamma)=2$.\qed

{\theorem $\Gamma(\mathbb{V})$ is complete if and only if $\mathbb{V}$ is one-dimensional.}\\
\pf Let $\Gamma(\mathbb{V})$ be complete. If possible, let $dim(\mathbb{V})>1$. Therefore, $\exists~ \alpha_1,\alpha_2 \in \mathbb{V}$ which is a basis or can be extended to a basis of $\mathbb{V}$. Then $\alpha_1$ and $\alpha_2$ are two non-adjacent vertices in $\Gamma(\mathbb{V})$, a contradiction. Therefore, $dim(\mathbb{V})=1$.

Conversely, let $\mathbb{V}$ be one-dimensional vector space generated by $\alpha$. Then any two non-null vectors $\mathbf{a}$ and $\mathbf{b}$ can be expressed as $c_1\alpha$ and $c_2\alpha$ respectively for non-zero $c_1,c_2 \in \mathbb{F}$ and hence $\mathbf{a}\sim \mathbf{b}$, thereby rendering the graph complete. \qed

{\theorem The domination number of $\Gamma(\mathbb{V}_\alpha)$ is 1.}\\
\pf The proof follows from the simple observation that $\alpha_1+\alpha_2+\cdots+\alpha_n$ is adjacent to all the vertices of $\Gamma(\mathbb{V}_\alpha)$. \qed

{\remark The set $\{\alpha_1,\alpha_2,\ldots,\alpha_n\}$ is a minimal dominating set of $\Gamma(\mathbb{V}_\alpha)$. Now, the question arises what is the maximum possible number of vertices in a minimal dominating set. The answer is given as $n$ in the next theorem.}

{\theorem If $D=\{\beta_1,\beta_2,\ldots,\beta_l\}$ is a minimal dominating set of $\Gamma(\mathbb{V}_\alpha)$, then $l \leq n$, i.e., the maximum cardinality of a minimal dominating set is $n$.}\\
\pf Since $D$ is a minimal dominating set, $\forall i \in \{1,2,\ldots,l\}, D_i=D \setminus \{\beta_i\}$ is not a dominating set. Therefore, $\forall i \in \{1,2,\ldots,l\}, \exists~ \gamma_i \in \Gamma(\mathbb{V}_\alpha)$ which is not adjacent to any element of $D_i$ but adjacent to $\beta_i$. Since, $\gamma_i \neq \theta, \exists~ \alpha_{t_i}$ such that $\gamma_i$ has non-zero component along $\alpha_{t_i}$. 

Now, as $\gamma_i$ is not adjacent to any element of $D_i$, so is $\alpha_{t_i}$. Thus, $\forall i \in \{1,2,\ldots,l\}, \exists~ \alpha_{t_i}$ such that $\alpha_{t_i} \sim \beta_i$, but $\alpha_{t_i} \not\sim \beta_k, \forall k\neq i$.

Claim: $i \neq j \Rightarrow \alpha_{t_i} \neq \alpha_{t_j}$.
Let, if possible, $i \neq j$ but $\alpha_{t_i} = \alpha_{t_j}$. As $\beta_i \sim \alpha_{t_i}$ and $\alpha_{t_i}=\alpha_{t_j}$, therefore $\beta_i \sim \alpha_{t_j}$. However, it contradicts that $\alpha_{t_i} \not\sim \beta_k, \forall k\neq i$. Hence, the claim.

As $\alpha_{t_1},\alpha_{t_2},\ldots,\alpha_{t_l}$ are all distinct, it follows that $l\leq n$. \qed

{\theorem \label{independence-number-theorem} The independence number of $\Gamma$ is $dim(\mathbb{V})$.}\\
\pf Since $\{\alpha_1,\alpha_2,\ldots,\alpha_n\}$ is an independent set in $\Gamma$, the independence number of $\Gamma \geq n=dim(\mathbb{V})$. Now, we show that any independent set can not have more than $n$ elements. Let, if possible, $\{\beta_1,\beta_2,\ldots,\beta_l\}$ be an independent set in $\Gamma$, where $l>n$. Since, $\beta_i \neq \theta, \forall i \in \{1,2,\ldots,l\}$, $\beta_i$ has at least one non-zero component along some $\alpha_{t_i}$, where $t_i \in \{1,2,\ldots,n\}$.\\
{\bf Claim:} $i \neq j \Rightarrow t_i \neq t_j$.\\
If $\exists i \neq j$ with $t_i=t_j=t$(say), then $\beta_i$ and $\beta_j$ has non-zero component along $\alpha_t$. This imply that $\beta_i \sim \beta_j$, a contradiction to the independence of $\beta_i$ and $\beta_j$. Hence, the claim is valid.

However, as there are exactly $n$ distinct $\alpha_i$'s, $l \leq n$, which is a contradiction. Thus, $\Gamma=n=dim(\mathbb{V})$.\qed

{\lemma \label{ind-imply-lin-ind} Let $I$ be an independent set in $\Gamma(\mathbb{V}_\alpha)$, then $I$ is a linearly independent subset of $\mathbb{V}$.}\\
\pf Let $I=\{\beta_1,\beta_2,\ldots,\beta_k\}$ be an independent set in $\Gamma$. By Theorem \ref{independence-number-theorem}, $k \leq n$. If possible, let $I$ be linearly dependent in $\mathbb{V}$. Then $\exists~ i \in \{1,2,\dots,k\}$ such that $\beta_i$ can be expressed as a linear combination of $\beta_1,\beta_2,\ldots,\beta_{i-1},\beta_{i+1},\ldots,\beta_k$, i.e.,
\begin{equation}\label{ind-set-equation}
\beta_i=c_1\beta_1+c_2\beta_2+\cdots+c_{i-1}\beta_{i-1}+c_{i+1}\beta_{i+1}+\cdots+
c_k\beta_k=\sum^{k}_{s=1,s\neq i}c_s\beta_s
\end{equation} 
Now, since $\{\alpha_1,\alpha_2,\ldots,\alpha_n\}$ is a basis of $\mathbb{V}$, let $\beta_j=\sum^{n}_{t=1} d_{tj}\alpha_t$ for $j=1,2,\ldots,i-1,i+1,\ldots k$. Thus, the expression of $\beta_i$ becomes
$$\beta_i=\sum^{k}_{s=1,s\neq i}c_s\sum^{n}_{t=1} d_{tj}\alpha_t=\sum^{n}_{t=1}d_t\alpha_t
\mbox{ for some scalars }d_t \in \mathbb{F}.$$
Since, $\beta_i \neq \theta$. Thus, $\beta_i$ has a non-zero component along some $\alpha_{t^*}$. 
Also, $\exists$ some $\beta_j, j\neq i$ such that $\beta_j$ has a non-zero component along $\alpha_{t^*}$. (as otherwise, if all $\beta_j,j\neq i$ has zero component along $\alpha_{t^*}$, then by Equation \ref{ind-set-equation}, $\beta_i$ has zero component along $\alpha_{t^*}$, which is not the case.)

Thus, $\beta_j \sim \beta_i$, a contradiction to the independence of $I$. Thus, $I$ is a linearly independent set in $\mathbb{V}$.\qed

{\remark Converse of Lemma \ref{ind-imply-lin-ind} is not true, in general. Consider a vector space $\mathbb{V}$, its basis $\{\alpha_1,\alpha_2,\alpha_3,\ldots,\alpha_n\}$ and the set $L=\{\alpha_1+\alpha_2,\alpha_2,\alpha_3,\ldots,\alpha_n\}$. Clearly $L$ is linearly independent in $\mathbb{V}$, but it is not an independent set in $\Gamma(\mathbb{V}_\alpha)$ as $\alpha_1+\alpha_2 \sim \alpha_2$.}

\section{Non-Zero Component Graph and Graph Isomorphisms}
In this section, we study the inter-relationship between the isomorphism of two vector spaces with the isomorphism of the two corresponding graphs. It is proved that two vector spaces are isomorphic if and only if their graphs are isomorphic. However, it is noted that a vector space isomorphism is a graph isomorphism (ignoring the null vector $\theta$), but a graph isomorphism may not be vector space isomorphism as shown in Example \ref{graph-iso-but-not-vecsp-iso}.

{\lemma \label{graph-isomorphic-implies-equal-dimension} Let $\mathbb{V}$ and $\mathbb{W}$ be two finite dimensional vector spaces over a field $\mathbb{F}$. If $\Gamma(\mathbb{V}_\alpha)$ and $\Gamma(\mathbb{W}_\beta)$ are isomorphic as graphs with respect to some basis $\{\alpha_1,\alpha_2,\ldots,\alpha_n\}$ and $\{\beta_1,\beta_2,\ldots,\beta_k\}$ of $\mathbb{V}$ and $\mathbb{W}$ respectively, then $dim(\mathbb{V})=dim(\mathbb{W})$, i.e., $n=k$.}\\
\pf Let $\varphi:\Gamma(\mathbb{V}_\alpha)\rightarrow \Gamma(\mathbb{W}_\beta)$ be a graph isomorphism. Since, $\alpha_1, \alpha_2,\ldots,\alpha_n$ is an independent set in $\Gamma(\mathbb{V}_\alpha)$, therefore $\varphi(\alpha_1),\varphi(\alpha_2),\ldots,\varphi(\alpha_n)$ is an independent set in $\Gamma(\mathbb{W}_\beta)$. Now, as in Theorem \ref{independence-number-theorem} it has been shown that cardinality of an independent set is less than or equal to the dimension of the vector space, it follows that $n \leq k$.

Again, $\varphi^{-1}:\Gamma(\mathbb{W}_\beta)\rightarrow \Gamma(\mathbb{V}_\alpha)$ is also a  graph isomorphism. Then, by similar arguments, it follows that $k \leq n$. Hence the lemma. \qed

{\theorem Let $\mathbb{V}$ and $\mathbb{W}$ be two finite dimensional vector spaces over a field $\mathbb{F}$. If $\mathbb{V}$ and $\mathbb{W}$ are isomorphic as vector spaces, then for any basis $\{\alpha_1,\alpha_2,\ldots,\alpha_n\}$ and $\{\beta_1,\beta_2,\ldots,\beta_n\}$ of $\mathbb{V}$ and $\mathbb{W}$ respectively, $\Gamma(\mathbb{V}_\alpha)$ and $\Gamma(\mathbb{W}_\beta)$ are isomorphic as graphs.}\\
\pf Let $\varphi: \mathbb{V}\rightarrow \mathbb{W}$ be a vector space isomorphism. Then $\{\varphi(\alpha_1),\varphi(\alpha_2),\ldots,\varphi(\alpha_n)\}$ is a basis of $\mathbb{W}$.
Consider the restriction $\overline{\varphi}$ of $\varphi$ on the non-null vectors of $\mathbb{V}$, i.e., $\overline{\varphi}: \Gamma(\mathbb{V}_\alpha) \rightarrow \Gamma(\mathbb{W}_{\varphi(\alpha)})$ given by $$\overline{\varphi}(a_1\alpha_1+a_2\alpha_2+\cdots+a_n\alpha_n)=a_1\varphi(\alpha_1)+a_2\varphi(\alpha_2)+\cdots
+a_n\varphi(\alpha_n)$$ where $a_i \in \mathbb{F}$ and $(a_1,a_2,\ldots,a_n)\neq(0,0,\ldots,0)$. Clearly, $\overline{\varphi}$ is a bijection. Now,\\ $\mathbf{a}\sim \mathbf{b} \mbox{ in }\Gamma(\mathbb{V}_\alpha) \Leftrightarrow \exists~ i \mbox{ such that }a_i \neq 0, b_i \neq 0 \Leftrightarrow \overline{\varphi}(\mathbf{a})\sim \overline{\varphi}(\mathbf{b})$.

Therefore, $\Gamma(\mathbb{V}_\alpha)$ and $\Gamma(\mathbb{W}_{\varphi(\alpha)})$ are graph isomorphic. Now, by Theorem \ref{basis-independent-theorem}, $\Gamma(\mathbb{W}_{\varphi(\alpha)})$ and $\Gamma(\mathbb{W}_{\beta})$ are isomorphic as graphs. Thus, $\Gamma(\mathbb{V}_\alpha)$ and $\Gamma(\mathbb{W}_{\beta})$ are isomorphic as graphs. \qed

{\theorem Let $\mathbb{V}$ and $\mathbb{W}$ be two finite dimensional vector spaces over a field $\mathbb{F}$. If for any basis $\{\alpha_1,\alpha_2,\ldots,$ $\alpha_n\}$ and $\{\beta_1,\beta_2,\ldots,\beta_k\}$ of $\mathbb{V}$ and $\mathbb{W}$ respectively, $\Gamma(\mathbb{V}_\alpha)$ and $\Gamma(\mathbb{W}_\beta)$ are isomorphic as graphs, then $\mathbb{V}$ and $\mathbb{W}$ are isomorphic as vector spaces.}\\
\pf Since $\Gamma(\mathbb{V}_\alpha)$ and $\Gamma(\mathbb{W}_\beta)$ are isomorphic as graphs, by Lemma \ref{graph-isomorphic-implies-equal-dimension}, $n=k$. Now, as $\mathbb{V}$ and $\mathbb{W}$ are finite dimensional vector spaces having same dimension over the same field $\mathbb{F}$, $\mathbb{V}$ and $\mathbb{W}$ are isomorphic as vector spaces. \qed  

{\example \label{graph-iso-but-not-vecsp-iso} Consider an one-dimensional vector space $\mathbb{V}$ over $\mathbb{Z}_5$ generated by $\alpha$ (say). Then $\Gamma(\mathbb{V}_\alpha)$ is a complete graph of 4 vertices with $V=\{\alpha,2\alpha,3\alpha,4\alpha\}$. Consider the map $T:\Gamma(\mathbb{V}_\alpha) \rightarrow \Gamma(\mathbb{V}_\alpha)$ given by $T(\alpha)=2\alpha,T(2\alpha)=\alpha,T(3\alpha)=4\alpha,T(4\alpha)=3\alpha$. Clearly, $T$ is a graph isomorphism, but as $T(2\alpha)=\alpha\neq 4\alpha=2(2\alpha)=2T(\alpha)$, $T$ is not linear.}

\section{Automorphisms of Non-Zero Component Graph}
In this section, we investigate the form of automorphisms of $\Gamma(\mathbb{V}_\alpha)$. It is shown that an automorphism maps $\{\alpha_1,\alpha_2,\ldots,\alpha_n\}$ to a basis of $\mathbb{V}$ of a special type, namely non-zero scalar multiples of a permutation of $\alpha_i$'s.

{\theorem \label{main-auto-theorem} Let $\varphi: \Gamma(\mathbb{V}_\alpha) \rightarrow \Gamma(\mathbb{V}_\alpha)$ be a graph automorphism. Then, $\varphi$ maps $\{\alpha_1,\alpha_2,\ldots,\alpha_n\}$ to another basis $\{\beta_1,\beta_2,\ldots,\beta_n\}$ such that there exists $\sigma \in S_n$, where each $\beta_i$ is of the form $c_i \alpha_{\sigma(i)}$ and each $c_i$'s are non-zero.}\\
\pf Let $\varphi: \Gamma(\mathbb{V}_\alpha) \rightarrow \Gamma(\mathbb{V}_\alpha)$ be a graph automorphism. Since, $\{\alpha_1,\alpha_2,\ldots,\alpha_n\}$ is an independent set of vertices in $\Gamma(\mathbb{V}_\alpha)$, therefore $\beta_i=\varphi(\alpha_i): i=1,2,\ldots,n$ is also an independent set of vertices in $\Gamma(\mathbb{V}_\alpha)$. Let $$\begin{array}{c}
\varphi(\alpha_1)=\beta_1=c_{11}\alpha_1+c_{12}\alpha_2+\cdots+c_{1n}\alpha_n\\
\varphi(\alpha_2)=\beta_2=c_{21}\alpha_1+c_{22}\alpha_2+\cdots+c_{2n}\alpha_n \\ \cdots \cdots \cdots \cdots \cdots \cdots \cdots \cdots \cdots \cdots \cdots \cdots \\
\varphi(\alpha_n)=\beta_n=c_{n1}\alpha_1+c_{n2}\alpha_2+\cdots+c_{nn}\alpha_n 
\end{array} $$
Since, $\beta_1 \neq \theta$ i.e., $\beta_1$ is not an isolated vertex, $\exists~ j_1 \in \{1,2,\ldots,n\}$ such that $c_{1j_1}\neq 0$. Therefore, $c_{ij_1}= 0, \forall i \neq 1$ (as $\beta_i$ is not adjacent to $\beta_1, \forall i \neq 1$.) Similarly, for $\beta_2$, $\exists~ j_2 \in \{1,2,\ldots,n\}$ such that $c_{2j_2}\neq 0$ and $c_{ij_2}= 0, \forall i \neq 2$. Moreover, $j_1 \neq j_2$ as $\beta_1$ and $\beta_2$ are not adjacent. Continuing in this manner, for $\beta_n$, $\exists~ j_n \in \{1,2,\ldots,n\}$ such that $c_{nj_n}\neq 0$ and $c_{ij_n}= 0, \forall i \neq n$ and $j_1,j_2,\ldots,j_n$ are all distinct numbers from $\{1,2,\ldots,n\}$.

Thus, $c_{kj_l}=0$ for $k \neq l$ and $c_{kj_k}\neq 0$, where $k,l \in \{1,2,\ldots,n\}$ and $j_1,j_2,\ldots,j_n$ is a permutation on $\{1,2,\ldots,n\}$. Set $\sigma=\left( \begin{array}{cccc}
1 & 2 & \cdots & n\\
j_1 & j_2 & \cdots & j_n
\end{array}\right)$. Therefore, $$\beta_i=c_{ij_i}\alpha_{j_i}=c_{ij_i}\alpha_{\sigma(i)},\mbox{ with }c_{ij_i}\neq 0,~~~~	 \forall i \in \{1,2,\ldots,n\}.$$ As $\{\alpha_1,\alpha_2,\ldots,\alpha_n\}$ is a basis, $\{\beta_1,\beta_2,\ldots,\beta_n\}$ is also a basis and hence the theorem. \qed

{\remark Although $\varphi$ maps the basis $\{\alpha_1,\alpha_2,\ldots,\alpha_n\}$ to another basis $\{\beta_1,\beta_2,\ldots,\beta_n\}$, it may not be a vector space isomorphism. It is because linearity of $\varphi$ is not guaranteed as shown in Example \ref{graph-iso-but-not-vecsp-iso}. However, the following result is true.}

{\theorem \label{special-auto-theorem} Let $\varphi$ be a graph automorphism, which maps $\alpha_i\mapsto c_{ij_i}\alpha_{\sigma(i)}$ for some $\sigma \in S_n$. Then, if $c \neq 0$, $\varphi(c\alpha_i)=d\alpha_{\sigma(i)}$ for some non-zero $d$. More generally, $\forall k \in \{1,2,\ldots,n\}$ if $c_1\cdot  c_2 \cdots c_k \neq 0$, then $$\varphi(c_1\alpha_{i_1}+c_2 \alpha_{i_2}+\cdots+c_k \alpha_{i_k})=d_1\alpha_{\sigma(i_1)}+d_2\alpha_{\sigma(i_2)}+\cdots+d_k \alpha_{\sigma(i_k)} $$ for some $d_i$'s with $d_1\cdot d_2\cdots d_k \neq 0$.}\\
\pf Since, $c\alpha_i \sim \alpha_i$, therefore $\varphi(c\alpha_i)\sim \varphi(\alpha_i)$ i.e., $\varphi(c\alpha_i)\sim c_{ij_i}\alpha_{\sigma(i)}$. Thus, $\varphi(c\alpha_i)$ has $\alpha_{\sigma(i)}$ as a non-zero component. If possible, let $\varphi(c\alpha_i)$ has a non-zero component along some other $\alpha_{\sigma(j)}$ for some $j\neq i$. Then $\varphi(c\alpha_i) \sim \alpha_{\sigma(j)}$ i.e., $\varphi(c\alpha_i) \sim \varphi(\alpha_j)$, which in turn implies $c\alpha_i \sim \alpha_j$ for $j\neq i$, a contradiction. Therefore, $\varphi(c\alpha_i)=d\alpha_{\sigma(i)}$ for some non-zero $d$.

For the general case, since
 $$c_1\alpha_{i_1}+c_2 \alpha_{i_2}+\cdots+c_k \alpha_{i_k} \sim \alpha_{i_1}$$
$$\Rightarrow \varphi(c_1\alpha_{i_1}+c_2 \alpha_{i_2}+\cdots+c_k \alpha_{i_k}) \sim \varphi(\alpha_{i_1})=c\alpha_{\sigma(i_1)} \mbox{ for some non-zero }c$$
$$\Rightarrow \varphi(c_1\alpha_{i_1}+c_2 \alpha_{i_2}+\cdots+c_k \alpha_{i_k})\mbox{ has a non-zero component along }\alpha_{\sigma(i_1)}$$
$$\Rightarrow \varphi(c_1\alpha_{i_1}+c_2 \alpha_{i_2}+\cdots+c_k \alpha_{i_k}) \sim \alpha_{\sigma(i_1)}~~~~~~~~~~~~~~~~~~~~~~~~~~~~~~~~~~~~~~~~$$
Similarly, 
$$\varphi(c_1\alpha_{i_1}+c_2 \alpha_{i_2}+\cdots+c_k \alpha_{i_k}) \sim \alpha_{\sigma(i_2)},\ldots,\varphi(c_1\alpha_{i_1}+c_2 \alpha_{i_2}+\cdots+c_k \alpha_{i_k}) \sim \alpha_{\sigma(i_k)}$$
Therefore, $\varphi(c_1\alpha_{i_1}+c_2 \alpha_{i_2}+\cdots+c_k \alpha_{i_k})=d_1\alpha_{\sigma(i_1)}+d_2\alpha_{\sigma(i_2)}+\cdots+d_k \alpha_{\sigma(i_k)} $ for some $d_i$'s with $d_1\cdot d_2\cdots d_k \neq 0$. \qed

{\corollary $\Gamma(\mathbb{V}_\alpha)$ is not vertex transitive if $dim(\mathbb{V})>1$.}\\
\pf Since $dim(\mathbb{V})\geq 2$, by Theorem \ref{special-auto-theorem}, there does not exist any automorphism which maps $\alpha_1$ to $\alpha_1+\alpha_2$. Hence, the result. \qed

\section{The Case of Finite Fields}
In this section, we find the degree of each vertices of $\Gamma(\mathbb{V})$ if the base field is finite. For more results, in the case of finite fields, please refer to \cite{angsu-jaa}.
{\remark \label{nbd-remark} The set of vertices adjacent to $\alpha_{i_1}+ \alpha_{i_2}+\cdots+ \alpha_{i_k}$ is same as the set of vertices adjacent to $c_1\alpha_{i_1}+c_2 \alpha_{i_2}+\cdots+c_k \alpha_{i_k}$ i.e., $N(\alpha_{i_1}+ \alpha_{i_2}+\cdots+ \alpha_{i_k})=N(c_1\alpha_{i_1}+c_2 \alpha_{i_2}+\cdots+c_k \alpha_{i_k})$ for $c_1 c_2 \cdots c_k \neq 0$.}

{\theorem Let $\mathbb{V}$ be a vector space over a finite field $\mathbb{F}$ with $q$ elements and $\Gamma$ be its associated graph with respect to a basis $\{\alpha_1,\alpha_2,\ldots,\alpha_n\}$. Then, the degree of the vertex $c_1\alpha_{i_1}+c_2 \alpha_{i_2}+\cdots+c_k \alpha_{i_k}$, where $c_1 c_2 \cdots c_k \neq 0$, is $(q^k -1)q^{n-k}-1$.}\\
\pf The number of vertices with $\alpha_{i_1}$ as non-zero component is $(q-1)q^{n-1}$ (including $\alpha_{i_1}$ itself). Therefore, $deg(\alpha_{i_1})=(q-1)q^{n-1}-1$.

The number of vertices with $\alpha_{i_1}$ or $\alpha_{i_2}$ as non-zero component is equal to
number of vertices with $\alpha_{i_1}$ as non-zero component $+$ number of vertices with $\alpha_{i_2}$ as non-zero component $-$ number of vertices with both $\alpha_{i_1}$ and $\alpha_{i_2}$ as non-zero component
$$=(q-1)q^{n-1}+(q-1)q^{n-1}-(q-1)^2 q^{n-2}=(q^2-1)q^{n-2}.$$
As this count includes the vertex $\alpha_{i_1}+\alpha_{i_2}$, $deg(\alpha_{i_1}+\alpha_{i_2})=(q^2-1)q^{n-2} - 1$.

Similarly, for finding the degree of $\alpha_{i_1}+\alpha_{i_2}+\alpha_{i_3}$, the number of vertices with $\alpha_{i_1}$ or $\alpha_{i_2}$ or $\alpha_{i_3}$ as non-zero component is equal to $$[(q-1)q^{n-1}+(q-1)q^{n-1}+(q-1)q^{n-1}]-[(q-1)^2 q^{n-2}+(q-1)^2 q^{n-2}+(q-1)^2 q^{n-2}]+(q-1)^3 q^{n-3}$$
$$=(q^3-1) q^{n-3}, \mbox{ and hence }deg(\alpha_{i_1}+\alpha_{i_2}+\alpha_{i_3})=(q^3-1) q^{n-3}-1.$$
Proceeding in this way, we get $$deg(\alpha_{i_1}+ \alpha_{i_2}+\cdots+ \alpha_{i_k})=(q^k -1)q^{n-k}-1.$$
Now, from Remark \ref{nbd-remark}, it follows that $$deg(c_1 \alpha_{i_1}+ c_2 \alpha_{i_2}+\cdots+c_k \alpha_{i_k})=(q^k -1)q^{n-k}-1.$$ \qed

\section{Conclusion}
In this paper, we represent a finite dimensional vector space as a graph and study various inter-relationships among $\Gamma(\mathbb{V})$ as a graph and $\mathbb{V}$ as a vector space. The main goal of these discussions was to study the nature of the automorphisms and establish the equivalence between the corresponding graph and vector space isomorphisms. Apart from this, we also study basic properties of completeness, connectedness, domination and independence number. As a topic of further research, one can look into the structure of maximal cliques and chromatic number of such graphs.
\section*{Acknowledgement}
The author is thankful to Bedanta Bose of University of Calcutta, Kolkata for bringing the manuscript in the final form. A special thanks goes to Dr. Usman Ali of Bahauddin Zakariya University, Pakistan for pointing out a mistake in an earlier version of the paper. The research is partially funded by NBHM Research Project Grant, (Sanction No. 2/48(10)/2013/ NBHM(R.P.)/R\&D II/695), Govt. of India.

\end{document}